\documentclass{article}[12pt]
\usepackage{bm, amscd,  amsfonts, amssymb, cite,texdraw}
\usepackage{amsmath, epsfig, cite, color, colordvi}
\usepackage{latexsym,graphicx}

\makeatletter \@addtoreset{equation}{section} \makeatother

\newtheorem{theo}{Theorem}[section]

\newtheorem{coro}[theo]{Corollary}

\def\k{\mathbf k}
\def\0{\mathbf 0}
\def\del{\boldsymbol \delta}
\def\S{{\mathbf S}}

\def\mid{{\,|\,}}

\def\pf{\noindent {\it Proof.\ }}
\def\qed{\hfill \rule{4pt}{7pt}}
\begin{document}

\begin{center}
{\Large The Method of Multiple Combinatorial Telescoping
% for two Identities \\[5pt]
%of Andrews on Upper Parity Indices of Partitions
}

\vskip 3mm
Daniel K. Du$^{1}$, Qing-Hu Hou$^{1,2}$ and Charles B. Mei$^2$
\vskip 3mm
$^1$Center for Applied Mathematics\\ Tianjin University,
Tianjin 300072, P. R. China\\[5pt]
$^2$Center for Combinatorics, LPMC-TJKLC \\
Nankai University, Tianjin 300071, P. R. China

\vskip 3mm

 E-mail: daniel@tju.edu.cn, hou@nankai.edu.cn, meib@mail.nankai.edu.cn

\end{center}

\begin{abstract}
We generalize the method of combinatorial telescoping to the case of multiple summations.
We shall demonstrate this idea by giving
combinatorial proofs for two identities of Andrews on parity indices of partitions.
\end{abstract}



\section{Introduction}\label{sec-intro}

The method of combinatorial telescoping for alternating sums was proposed by Chen et al.\cite{Chen2011},
which can be used to show that an alternating sum satisfies certain recurrence relation combinatorially. With this method, one  give combinatorial interpretations for many $q$-series identities such as Watson's identity\cite{Watson} and Sylvester's identity\cite{Syl}.
For $q$-series identities on positive terms, Chen et al.\cite{Chen2012} presented the
corresponding combinatorial telescoping, based on which they established a
combinatorial proof for an identity due to Andrews \cite{Andrews2010}.

In this paper, we shall generalize the method of combinatorial telescoping to the multiple cases.
More precisely, we shall give the combinatorial telescoping for $q$-series identities of the following form
\begin{equation}\label{iden-multiple}
\sum_{\k=\0}^\infty (-1)^{\del \cdot \k} f(\k) =
\sum_{\k=\0}^\infty (-1)^{\del \cdot \k} g(\k),
\end{equation}
where $\k=(k_1,\ldots,k_m)$ and $\del=(\delta_1,\ldots,\delta_m)$
are $m$-dimensional vectors, and $\del_i \in\{0,1\}$.

Assume that
$f(\k)$ and $g(\k)$ are weighted counts of sets
$A_{\k}$ and $B_{\k}$, respectively, that is,
\[
  f(\k) = \sum_{\alpha\in A_{\k}} w(\alpha) \quad \mbox{ and } \quad
  g(\k) = \sum_{\alpha\in B_{\k}}w(\alpha).
\]
Motivated by the idea of creative telescoping of Zeilberger\cite{Knuth1994,A=B1996,Zeilberger1991},
we will construct sets
 $\{H_{i,\k}\}_{i=1}^m$ with a weight assignment $w$ such that there exists
a sequence of weight preserving bijections
\begin{equation}\label{biject-multiple}
 \phi_{\k} :
     A_{\k} \bigcup_{\{i| \del_i=0\}}  H_{i,\S_i \k}
     \longrightarrow   B_{\k}
      \bigcup_{i=1}^m H_{i,\k} \bigcup_{\{i| \del_i=1\}} H_{i,\S_i \k},
\end{equation}
where $\S_i$ is the shift operator on the $i$-th part, i.e.,
\[
 \S_i \k = ( k_1,\ldots,k_{i-1},k_i+1,k_{i+1},\ldots,k_m).
\]

Since $\phi _{\k}$ and $\phi_{\S_i \k}$ are weight preserving for
$1\le i\le m$, both
\[
 \phi_{\k}(H_{i,\S_i \k}) \quad \mbox{ and } \quad
 \phi_{\S_i\k}(H_{i,\S_i \k})
\]
have the same weight as $H_{i,\S_i \k}$. Realizing that for all $1\le i \le m$ with $\del_i=1$, we have
\[
 \phi_{\k}^{-1}(H_{i,\S_i \k}) \subseteq
 A_{\k} \bigcup_{\{i| \del_i=0\}}  H_{i,\S_i \k}
\]
and
\[
 \phi_{\S_i \k}^{-1}(H_{i,\S_i \k}) \subseteq
 A_{\S_i\k} \bigcup_{\{i| \del_i=0\}}  H_{i,\S_i^2 \k},
\]
which implies that the corresponding summands will cancel each other
in the desired identity \eqref{iden-multiple}. More precisely, if we set
\[
 h_i(\k) = \sum_{\alpha\in H_{i,\k}} w (\alpha), \quad (1\le i\le m),
\]
then the bijection \eqref{biject-multiple} implies that
\begin{eqnarray}\label{rec-multiple}
 &  & \nonumber f(\k)+ \sum_{\{i| \del_i=0 \}} h_i(\S_i \k)\\
 &=& \nonumber g(\k)+ \sum_{i=1}^m h_i(\k) +
             \sum_{\{i| \del_i=1 \}} h_i(\S_i \k)\\
 &=& g(\k)+ \sum_{\{i| \del_i=0 \}}h_i(\k) +
         \sum_{\{i| \del_i=1 \}}h_i(\k) + \sum_{\{i| \del_i=1 \}} h_i(\S_i \k).
\end{eqnarray}
We assume, like the conditions for the creative telescoping\cite{Knuth1994,A=B1996,Zeilberger1991},
that $H_i({\bold 0})=\emptyset$ and $H_i(\k)$ vanishes for sufficiently
large $\k$ for $1\le i \le m$. Multiplying $(-1)^{\del\cdot\k}$  and summing  over
 $\k$ on both sides of \eqref{rec-multiple}, since 
\[
(-1)^{\del\cdot\k} h_i(\S_i\k) + (-1)^{\del\cdot\S_i\k} h_i(\S_i\k) = 0
\]
for $1\le i \le m$ with $\del_i=1$, we will obtain the identity \eqref{iden-multiple}, which
is often an identity we wish to establish.

Indeed, once we have bijections $\phi_{\k}$ in \eqref{biject-multiple}, 
combining all these bijections, we are lead to a correspondence
\begin{equation}\label{eq-phi}
\phi: A\cup C \rightarrow B\cup C,
\end{equation}
where
\[
  A= \bigcup_{\k=0}^\infty A_{\k},\quad
 B = \bigcup_{\k=0}^\infty B_{\k}, \quad \mbox{ and }\quad
 C = \bigcup_{\{i|\del_i=0\}} H_{i,\S_i \k}.
\]
To be more specific, we can derive a bijection
\[
 \phi : A\cup C\setminus  D
   \longrightarrow   B\cup C,
\]
where
\[
 D =\bigcup_{\k=\0}^\infty \bigcup_{\{i| \del_i=1\}}\phi_{\k}^{-1}\left(
     (H_{i, \k}\cup H_{i,\S_i \k})\right),
\]
and an involution
\[
 \psi : D\longrightarrow D
\]
given by $\phi(\alpha)=\phi_{\k}(\alpha)$ if $\alpha \in A_{\k} \cup
C$ and
\[
 \psi(\alpha) = \left\{
\begin{array}{ll}
 \phi^{-1}_{\S_i^{-1} \k}\phi_{\k}(\alpha), &
    \mbox{ if } \alpha \in \phi^{-1}(H_{i,\k}),\\[5pt]
 \phi^{-1}_{\S_i \k}\phi_{\k}(\alpha), &
    \mbox{ if } \alpha \in \phi^{-1}(H_{i, \S_i \k}).\\
\end{array}
\right.
\]
As shown in \cite{Chen2012}, using the method of cancelation (see \cite{Feldman1995}), the
bijection $\phi$ in \eqref{eq-phi} implies a bijection
\[
\psi: A \rightarrow B,
\]
which gives a combinatorial interpretation of
the desired $q$-series identity \eqref{iden-multiple}.

The above approach to proving an identity like \eqref{iden-multiple} is called
\emph{multiple combinatorial telescoping}. To illustrate the idea of this method,
we shall prove two
$q$-series identiies proposed by Andrews\cite{Andrews2010}.

In the study of parity in partition identities, Andrews\cite{Andrews2010} proposed
fifteen problems. Two of them, labeled as Question $9$ and Question $10$ in \cite{Andrews2010},
asked for proving the following two identities of sum on double variables:
\begin{align}
\sum_{n,k\geq
0}\frac{(-1)^n
q^{(n-k)^2+k^2+n-k}}{(-q;q)_n(q;q)_{2k-1}(q;q)_{n-2k+1}} = \sum_{n=1}^{\infty}(-1)^nq^{n^2}, \label{eq-yee-9}\\[9pt]
\sum_{n,k\geq 0}\frac{(-1)^n
q^{(n-k)^2+k^2+n+k}}{(-q;q)_n(q;q)_{2k}(q;q)_{n-2k}} =\sum_{n=0}^{\infty}(-1)^nq^{n^2}. \label{eq-yee-10}
\end{align}
Recently Yee \cite{Yee2010} and Chu\cite{Chu2012} provided algebraic proofs for \eqref{eq-yee-9} and \eqref{eq-yee-10} independently, while
the combinatorial interpretation is still open.

In the framework of the method of multiple combinatorial telescoping,
we shall give a more extensive result as follows.
\begin{theo}\label{the-main}
\begin{eqnarray}
\sum_{m,k\geq 0}\frac{(-a)^m
q^{(m-k)^2+k^2+m-k}}{(aq^2;q^2)_m}\left[\!
\begin{array}{c}
  m \\
  2k-1
\end{array}\!
\right]_q = \sum_{n\geq 1} (-a)^n q^{n^2},\label{eq-q9a}\\[5pt]
\sum_{m,k\geq 0}\frac{(-a)^m
q^{(m-k)^2+k^2+k}}{(aq^2;q^2)_m}\left[\!
\begin{array}{c}
  m \\
  2k
\end{array}\!
\right]_q = \sum_{n\geq 0} (-a)^n q^{n^2},\label{eq-q10a}
\end{eqnarray}
\end{theo}
By setting $a=1$ in \eqref{eq-q9a} and \eqref{eq-q10a},
they are reduce to \eqref{eq-yee-9} and \eqref{eq-yee-10}, respectively,
which give combinatorial answers to the open questions of Andrews.

\section{Multiple Combinatorial Telescopings for Identities of Andrews}\label{sec-proof9}
In this section, by constructing multiple combinatorial telescopings, we shall
give certain bijections for recurrence relations of identity \eqref{eq-q9a} and \eqref{eq-q10a},
respectively. Based on the method of cancelation (see \cite{Feldman1995}), these bijections
lead to corresponding combinatorial interpretations of Question $9$ and Question $10$
proposed by Andrews in \cite{Andrews2010}. Since the proofs of \eqref{eq-q9a} and \eqref{eq-q10a}
are similar, we introduce the procedure for \eqref{eq-q9a} and a sketch proof for \eqref{eq-q10a}.

Let us recall some definitions concerning partitions as used in Andrews\cite{Andrews1998}.
A partition is
a non-increasing finite sequence of positive integers
\[
\lambda=(\lambda_1,\ldots,\lambda_l).
\]
 The integers $\lambda_i$ are
called the parts of $\lambda$. The sum of parts and the number of
parts are denoted by
\[
\mid \lambda \mid =\lambda_1+ \dots
+\lambda_l,
\]
and $\ell(\lambda)=l$, respectively. The number of
k-parts in $\lambda$ is denoted by $m_k(\lambda)$. The special
partition with no parts is denoted by $\emptyset$. We shall use
diagrams to represent partitions and use rows to represent parts.

\begin{figure}[h,t]
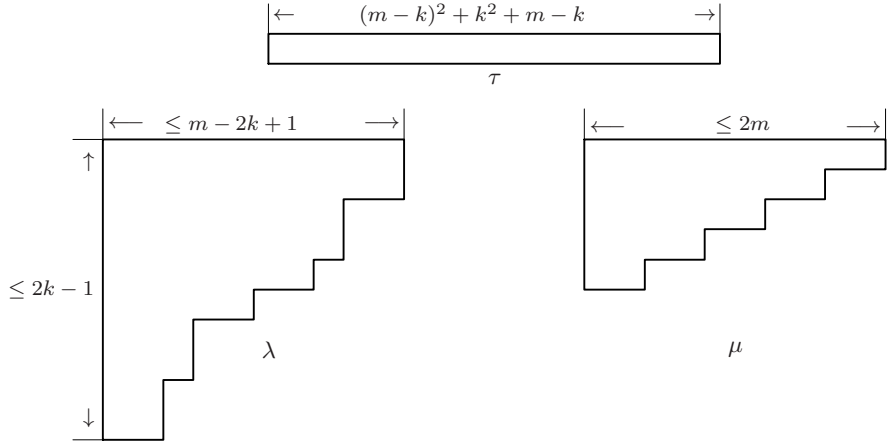

\centertexdraw{ \drawdim mm \linewd 0.25 \setgray 0

% \tau
\move(30 50)\lvec(90 50)\lvec(90 54)\lvec(30 54)\lvec(30 50)

\textref h:C v:C \htext(57 56.5){\footnotesize$(m-k)^2+k^2+m-k$}
 \textref h:C v:C \htext(32 56.5){\footnotesize$\leftarrow$}
  \textref h:C v:C \htext(88 56.5){\footnotesize$\rightarrow$}

% \lambda
\move(8 0)\lvec(16 0)\lvec(16 8)\lvec(20 8)\lvec(20 16)\lvec(28
16)\lvec(28 20)\lvec(36 20)\lvec(36 24)\lvec(40 24)\lvec(40 32)
\lvec(48 32)\lvec(48 40)\lvec(8 40)\lvec(8 0)

 \textref h:C v:C \htext(25 42){\footnotesize$\leq m-2k+1$}
 \textref h:C v:C \htext(11 42){\footnotesize$\longleftarrow$}
 \textref h:C v:C \htext(45 42){\footnotesize$\longrightarrow$}

  \textref h:C v:C \htext(1 20){\footnotesize$\leq 2k-1$}
 \textref h:C v:C \htext(6 37){\footnotesize$\uparrow$}
 \textref h:C v:C \htext(6 3){\footnotesize$\downarrow$}

%\mu
\move(72 40)\lvec(112 40)\lvec(112 36)\lvec(104 36)\lvec(104
32)\lvec(96 32)\lvec(96 28)\lvec(88 28)\lvec(88 24)\lvec(80 24)\lvec(80 20)
\lvec(72 20)\lvec(72 40)

 \textref h:C v:C \htext(93 42){\footnotesize$\leq 2m$}
 \textref h:C v:C \htext(75 41.5){\footnotesize$\longleftarrow$}
  \textref h:C v:C \htext(109 41.5){\footnotesize$\longrightarrow$}

\linewd 0.1 \setgray 0
%tau
\move(30 54)\lvec(30 58)
\move(90 54)\lvec(90 58)
%lambda
\move(4 0)\lvec(8 0)
\move(4 40)\lvec(8 40)
\move(8 40)\lvec(8 44)
\move(48 40)\lvec(48 44)
%mu
\move(112 44)\lvec(112 40)
\move(72 44)\lvec(72 40)

\textref h:C v:C \htext(60 48){\small$\tau$}
 \textref h:C v:C \htext(30 12){\small$\lambda$}
  \textref h:C v:C \htext(92 12){\small$\mu$} }
  \caption{Illustration of an element $(\tau,\lambda,\mu)\in P_{m,k}$.}
  \label{fig-9-case0}
\end{figure}

In order to prove identity \eqref{eq-q9a}, we let
\[
F(q) = \sum_{m,k}\frac{(-a)^m
q^{(m-k)^2+k^2+m-k}}{(aq^2;q^2)_m}\left[\!
\begin{array}{c}
  m \\
  2k-1
\end{array}\!
\right]_q,
\]
and the corresponding summand
\begin{equation}\label{eq-eq9-summation}
F_{m,k} = \frac{(-a)^m
q^{(m-k)^2+k^2+m-k}}{(aq^2;q^2)_m}\left[\!
\begin{array}{c}
  m \\
  2k-1
\end{array}\!
\right]_q.
\end{equation}
Define $P_{m,k}$ to be the set of triples $(\tau,\lambda,\mu)$,
where $\tau$ is a partition with only one part
$(m-k)^2+k^2+m-k$, $\lambda$ is a partition with no more than $2k-1$ parts and each part not exceeding $m-2k+1$, 
and $\mu$ is a partition with only even not exceeding $2m$; see Figure \ref{fig-9-case0}.
In particular, when $m=k=0$, we have
$P_{0,0}=\emptyset$.  Moreover,
one can see that the $(m,k)$-th summand $F_{m,k}$ of $F(q)$ in \eqref{eq-eq9-summation} can be
viewed as the weight of $P_{m,k}$, that is,
\[
\sum_{(\tau,\lambda,\mu)\in
P_{m,k}}a^{m+\ell(\mu)}q^{\mid \tau \mid +\mid \lambda \mid+\mid \mu \mid}.
\]
According to the exponent of $a$ in the above definition, we divide
$P_{m,k}$ into a disjoint union of subsets
$$P_{n,m,k}=\{(\tau,\lambda,\mu )\in P_{m,k}: \ell(\mu)=n-m\},$$
with
$P_{n,m,k}=\emptyset$ when $m=k=0$ or $2k-1>m>n$ .
By constructing three sets of triples $(\tau,\lambda,\mu )$ as follows
\[
\begin{array}{l}
G_{n,m,k}=\{(\tau,\lambda,\mu )\in
P_{n,m,k}:\ell(\lambda)=2k-1,m_2(\mu)=0\},\\[5pt]
H_{n,m,k}=\{(\tau,\lambda,\mu )\in
P_{n,m,k}:\ell(\lambda)=2k-2,m_{2m}(\mu)=0\},\\[5pt]
K_{n,m,k}=\{(\tau,\lambda,\mu )\in P_{n,m,k}:\ell(\lambda)\leq
2k-3,m_{2m}(\mu)=0\},
\end{array}
\]
we have the following combinatorial telescoping relation for $P_{n,m,k}$.
\begin{theo}\label{thm-mct}
For any positive integer $n$ and nonnegative integers $m$ and $k$, there is a bijection
\[
\begin{array}{l}
 \phi_{n,m,k}:  \\[5pt]
\quad P_{n,m,k}\cup \{2n-1\} \times
P_{n-1,m,k}\cup K_{n,m+1,k} \rightarrow \\[5pt]
  \quad G_{n,m,k}\cup
G_{n,m+1,k}\cup  H_{n,m,k}\cup H_{n,m+1,k}\cup
  K_{n,m,k}\cup K_{n,m+1,k+1}  \cup K_{n,m+1,k}.
\label{rec-q9}
\end{array}
\]
\end{theo}

\pf Let
\[
U_{n,m,k}=\{(\tau,\lambda,\mu )\in
P_{n,m,k}:\ell(\lambda)\leq 2k-2,m_{2m}(\mu)\neq0\},
\]
and
\[
T_{n,m,k}=\{(\tau,\lambda,\mu )\in
P_{n,m,k}:\ell(\lambda)=2k-1,m_2(\mu) \neq 0\},
\]
are two sets of triples $(\tau,\lambda,\mu )$, then we can divide $P_{n,m,k}$ into 
a disjoint union of five subsets, that is
\[
P_{n,m,k}=G_{n,m,k} \cup H_{n,m,k}\cup K_{n,m,k} \cup U_{n,m,k} \cup
T_{n,m,k}.
\]
Now we construct bijections $\phi_{n,m,k}$ like \eqref{biject-multiple}. To this purpose, we can classify
\[
P_{n,m,k}\cup \{2n-1\} \times P_{n-1,m,k}\cup K_{n,m+1,k}
\]
into four cases as below. We shall show that the elements in the first case are fix points, while the bijections for other three cases as follows:
\begin{equation*}
\begin{array}{lrll}
\varphi_1: &\{2n-1\}\times F_{n-1,m,k}&\rightarrow & G_{n,m+1,k}  ,\\
\varphi_2: &U_{n,m,k} &\rightarrow &  H_{n,m+1,k} ,\\
\varphi_3: &T_{n,m,k} &\rightarrow & K_{n,m+1,k+1} .
\end{array}
\end{equation*}
To be more specific, we have the following four cases.

Case $0$. For $(\tau,\lambda,\mu )\in G_{n,m,k}\cup H_{n,m,k} \cup
K_{n,m,k}\cup K_{n,m+1,k}$, let the image be itself.

\begin{figure}[h,t]
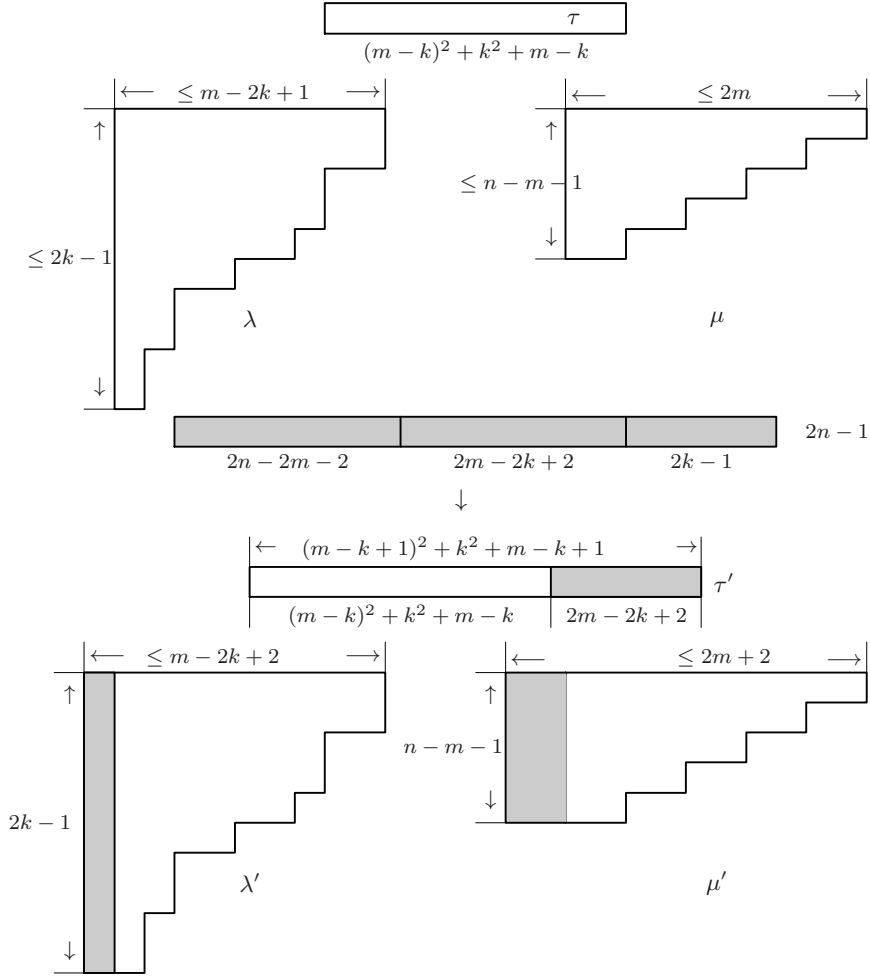

\centertexdraw{ \drawdim mm \linewd 0.25 \setgray 0
% \tau
\move(40 125)\lvec(80 125)\lvec(80 129)\lvec(40 129)\lvec(40 125)

\textref h:C v:C \htext(60 122.5){\footnotesize$(m-k)^2+k^2+m-k$}

% \lambda
\move(12 75)\lvec(16 75)\lvec(16 83)\lvec(20 83)\lvec(20 91)\lvec(28
91)\lvec(28 95)\lvec(36 95)\lvec(36 99)\lvec(40 99)\lvec(40 107)
\lvec(48 107)\lvec(48 115)\lvec(12 115)\lvec(12 75)

 \textref h:C v:C \htext(29 117){\footnotesize$\leq m-2k+1$}
 \textref h:C v:C \htext(15 117){\footnotesize$\longleftarrow$}
 \textref h:C v:C \htext(45 117){\footnotesize$\longrightarrow$}

  \textref h:C v:C \htext(6 95){\footnotesize$ \leq 2k-1$}
 \textref h:C v:C \htext(10 112){\footnotesize$\uparrow$}
 \textref h:C v:C \htext(10 78){\footnotesize$\downarrow$}

%\mu
\move(72 115)\lvec(112 115)\lvec(112 111)\lvec(104 111)\lvec(104
107)\lvec(96 107)\lvec(96 103)\lvec(88 103)\lvec(88 99)\lvec(80 99)\lvec(80 95)
\lvec(72 95)\lvec(72 115)

 \textref h:C v:C \htext(93 117){\footnotesize$\leq 2m$}
 \textref h:C v:C \htext(75 116.5){\footnotesize$\longleftarrow$}
  \textref h:C v:C \htext(109 116.5){\footnotesize$\longrightarrow$}

   \textref h:C v:C \htext(66 105){\footnotesize$ \leq n-m-1$}
 \textref h:C v:C \htext(70 112){\footnotesize$\uparrow$}
  \textref h:C v:C \htext(70 98){\footnotesize$\downarrow$}

%2n-1
\move(20 70)\lvec(100 70)\lvec(100 74)\lvec(20 74)\lvec(20 70)\ifill f:0.8
\move(20 70)\lvec(100 70)\lvec(100 74)\lvec(20 74)\lvec(20 70)
\move(50 70)\lvec(50 74)
\move(80 70)\lvec(80 74)
\textref h:C v:C \htext(35 68){\footnotesize$2n-2m-2$}
\textref h:C v:C \htext(65 68){\footnotesize$2m-2k+2$}
\textref h:C v:C \htext(90 68){\footnotesize$2k-1$}
\textref h:C v:C \htext(108 72){\footnotesize$2n-1$}

%downarrow
\textref h:C v:C \htext(58 63){\small$\downarrow$}

% \tau'
\move(30 50)\lvec(90 50)\lvec(90 54)\lvec(30 54)\lvec(30 50)
\move(70 50)\lvec(90 50)\lvec(90 54)\lvec(70 54)\lvec(70 50)\ifill f:0.8
\move(70 50)\lvec(90 50)\lvec(90 54)\lvec(70 54)\lvec(70 50)

\textref h:C v:C \htext(57 56.5){\footnotesize$(m-k+1)^2+k^2+m-k+1$}
 \textref h:C v:C \htext(32 56.5){\footnotesize$\leftarrow$}
  \textref h:C v:C \htext(88 56.5){\footnotesize$\rightarrow$}

\textref h:C v:C \htext(50 47.5){\footnotesize$(m-k)^2+k^2+m-k$}
%\textref h:C v:C \htext(57 44.5){\footnotesize$(m-k+1)^2+k^2+m-k+1$}
\textref h:C v:C \htext(80 47.5){\footnotesize$2m-2k+2$}

% \lambda'
\move(8 0)\lvec(16 0)\lvec(16 8)\lvec(20 8)\lvec(20 16)\lvec(28
16)\lvec(28 20)\lvec(36 20)\lvec(36 24)\lvec(40 24)\lvec(40 32)
\lvec(48 32)\lvec(48 40)\lvec(8 40)\lvec(8 0)
\move(8 0)\lvec(12 0)\lvec(12 40)\lvec(8 40)\lvec(8 0)\ifill f:0.8
\move(8 0)\lvec(12 0)\lvec(12 40)\lvec(8 40)\lvec(8 0)

 \textref h:C v:C \htext(25 42){\footnotesize$\leq m-2k+2$}
 \textref h:C v:C \htext(11 42){\footnotesize$\longleftarrow$}
 \textref h:C v:C \htext(45 42){\footnotesize$\longrightarrow$}

  \textref h:C v:C \htext(2 20){\footnotesize$ 2k-1$}
 \textref h:C v:C \htext(6 37){\footnotesize$\uparrow$}
 \textref h:C v:C \htext(6 3){\footnotesize$\downarrow$}

%\mu'
\move(72 40)\lvec(112 40)\lvec(112 36)\lvec(104 36)\lvec(104
32)\lvec(96 32)\lvec(96 28)\lvec(88 28)\lvec(88 24)\lvec(80 24)\lvec(80 20)
\lvec(72 20)\lvec(72 40)
\move(72 40)\lvec(64 40)\lvec(64 20)\lvec(72 20)\ifill f:0.8
\move(72 40)\lvec(64 40)\lvec(64 20)\lvec(72 20)

 \textref h:C v:C \htext(93 42){\footnotesize$\leq 2m+2$}
 \textref h:C v:C \htext(67 41.5){\footnotesize$\longleftarrow$}
  \textref h:C v:C \htext(109 41.5){\footnotesize$\longrightarrow$}

   \textref h:C v:C \htext(57 30){\footnotesize$ n-m-1$}
 \textref h:C v:C \htext(62 37){\footnotesize$\uparrow$}
  \textref h:C v:C \htext(62 23){\footnotesize$\downarrow$}

\linewd 0.1 \setgray 0
%lambda
\move(12 75)\lvec(8 75)
\move(12 115)\lvec(8 115)
\move(12 115)\lvec(12 119)
\move(48 115)\lvec(48 119)
%mu
\move(72 95)\lvec(68 95)
\move(72 115)\lvec(68 115)
\move(72 115)\lvec(72 119)
\move(112 115)\lvec(112 119)
%tau'
\move(30 54)\lvec(30 58)
\move(90 54)\lvec(90 58)
\move(30 50)\lvec(30 46)
\move(90 50)\lvec(90 46)
\move(70 50)\lvec(70 46)
%lambda'
\move(4 0)\lvec(8 0)
\move(4 40)\lvec(8 40)
\move(8 40)\lvec(8 44)
\move(48 40)\lvec(48 44)
%mu'
\move(112 44)\lvec(112 40)
\move(64 44)\lvec(64 40)
\move(64 40)\lvec(60 40)
\move(64 20)\lvec(60 20)

\textref h:C v:C \htext(73 127){\small$\tau$}
 \textref h:C v:C \htext(30 87){\small$\lambda$}
  \textref h:C v:C \htext(92 87){\small$\mu$}

\textref h:C v:C \htext(93 52){\small$\tau'$}
 \textref h:C v:C \htext(30 12){\small$\lambda'$}
  \textref h:C v:C \htext(92 12){\small$\mu'$} }

  \caption{The bijection $\varphi_1$ in Case $1$.}\label{fig-9-case1}
\end{figure}

 Case $1$. For $(\tau,\lambda,\mu )\in P_{n-1,m,k}$, as $\ell(\lambda)=2k-1$, we add
 a column equals $2k-1$
 to  $\lambda$ and  obtain a new partition
 $\lambda'$ ,where $\ell(\lambda')= 2k-1$ and $\lambda'_1\leq
 m-2k+2$. By adding $2m-2k+2$ to $\tau$ we get
 $\tau'=((m-k+1)^2+k^2+m-k+1)$. And we can add 2 columns equal $n-m-1$ to $\mu$
 without changing the number of parts , thus the new
 partition $\mu'$ satisfies that $\ell(\mu')=n-m-1$, $\mu'_1 \leq 2m+2$ and $m_2(\mu')=0$.
 We see that the weight of $(\tau,\lambda,\mu )$ is less than $(\tau',\lambda',\mu'
 )$ by $2n-1$. So we obtain the bijection $\varphi_1:  \{2n-1\}\times P_{n-1,m,k} \rightarrow G_{n,m+1,k}$
 defined by $(2n-1,(\tau,\lambda,\mu )) \mapsto (\tau',\lambda',\mu'
 )$. Figure \ref{fig-9-case1} gives an illustration of the correspondence $\varphi_1$.

\begin{figure}[h,t]
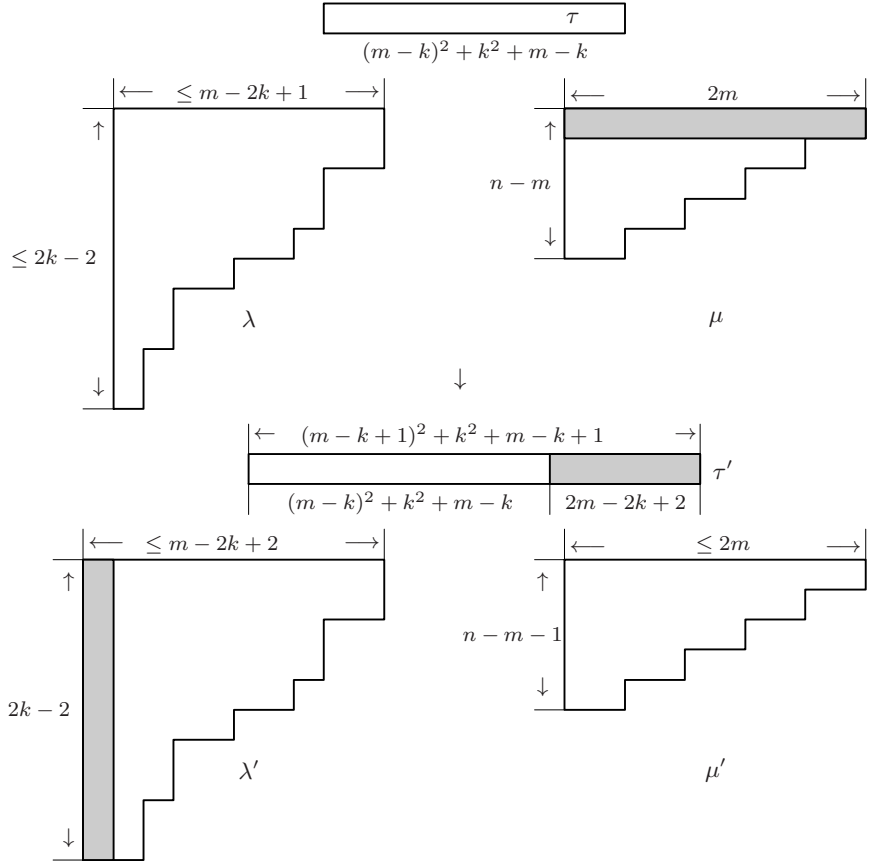

\centertexdraw{ \drawdim mm \linewd 0.25 \setgray 0
% \tau
\move(40 110)\lvec(80 110)\lvec(80 114)\lvec(40 114)\lvec(40 110)

\textref h:C v:C \htext(60 107.5){\footnotesize$(m-k)^2+k^2+m-k$}

% \lambda
\move(12 60)\lvec(16 60)\lvec(16 68)\lvec(20 68)\lvec(20 76)\lvec(28
76)\lvec(28 80)\lvec(36 80)\lvec(36 84)\lvec(40 84)\lvec(40 92)
\lvec(48 92)\lvec(48 100)\lvec(12 100)\lvec(12 60)

 \textref h:C v:C \htext(29 102){\footnotesize$\leq m-2k+1$}
 \textref h:C v:C \htext(15 102){\footnotesize$\longleftarrow$}
 \textref h:C v:C \htext(45 102){\footnotesize$\longrightarrow$}

  \textref h:C v:C \htext(4 80){\footnotesize$ \leq 2k-2$}
 \textref h:C v:C \htext(10 97){\footnotesize$\uparrow$}
 \textref h:C v:C \htext(10 63){\footnotesize$\downarrow$}

%\mu
\move(72 100)\lvec(112 100)\lvec(112 96)\lvec(104 96)\lvec(104
92)\lvec(96 92)\lvec(96 88)\lvec(88 88)\lvec(88 84)\lvec(80 84)\lvec(80 80)
\lvec(72 80)\lvec(72 100)
\move(72 100)\lvec(112 100)\lvec(112 96)\lvec(72 96)\lvec(72 100)\ifill f:0.8
\move(72 100)\lvec(112 100)\lvec(112 96)\lvec(72 96)\lvec(72 100)

 \textref h:C v:C \htext(93 102){\footnotesize$ 2m$}
 \textref h:C v:C \htext(75 101.5){\footnotesize$\longleftarrow$}
  \textref h:C v:C \htext(109 101.5){\footnotesize$\longrightarrow$}

   \textref h:C v:C \htext(66 90){\footnotesize$ n-m$}
 \textref h:C v:C \htext(70 97){\footnotesize$\uparrow$}
  \textref h:C v:C \htext(70 83){\footnotesize$\downarrow$}

%downarrow
\textref h:C v:C \htext(58 64){\small$\downarrow$}

% \tau'
\move(30 50)\lvec(90 50)\lvec(90 54)\lvec(30 54)\lvec(30 50)
\move(70 50)\lvec(90 50)\lvec(90 54)\lvec(70 54)\lvec(70 50)\ifill f:0.8
\move(70 50)\lvec(90 50)\lvec(90 54)\lvec(70 54)\lvec(70 50)

\textref h:C v:C \htext(57 56.5){\footnotesize$(m-k+1)^2+k^2+m-k+1$}
 \textref h:C v:C \htext(32 56.5){\footnotesize$\leftarrow$}
  \textref h:C v:C \htext(88 56.5){\footnotesize$\rightarrow$}

\textref h:C v:C \htext(50 47.5){\footnotesize$(m-k)^2+k^2+m-k$}
%\textref h:C v:C \htext(57 44.5){\footnotesize$(m-k+1)^2+k^2+m-k+1$}
\textref h:C v:C \htext(80 47.5){\footnotesize$2m-2k+2$}

% \lambda'
\move(8 0)\lvec(16 0)\lvec(16 8)\lvec(20 8)\lvec(20 16)\lvec(28
16)\lvec(28 20)\lvec(36 20)\lvec(36 24)\lvec(40 24)\lvec(40 32)
\lvec(48 32)\lvec(48 40)\lvec(8 40)\lvec(8 0)
\move(8 0)\lvec(12 0)\lvec(12 40)\lvec(8 40)\lvec(8 0)\ifill f:0.8
\move(8 0)\lvec(12 0)\lvec(12 40)\lvec(8 40)\lvec(8 0)

 \textref h:C v:C \htext(25 42){\footnotesize$\leq m-2k+2$}
 \textref h:C v:C \htext(11 42){\footnotesize$\longleftarrow$}
 \textref h:C v:C \htext(45 42){\footnotesize$\longrightarrow$}

  \textref h:C v:C \htext(2 20){\footnotesize$ 2k-2$}
 \textref h:C v:C \htext(6 37){\footnotesize$\uparrow$}
 \textref h:C v:C \htext(6 3){\footnotesize$\downarrow$}

%\mu'
\move(72 40)\lvec(112 40)\lvec(112 36)\lvec(104 36)\lvec(104
32)\lvec(96 32)\lvec(96 28)\lvec(88 28)\lvec(88 24)\lvec(80 24)\lvec(80 20)
\lvec(72 20)\lvec(72 40)

 \textref h:C v:C \htext(93 42){\footnotesize$\leq 2m$}
 \textref h:C v:C \htext(75 41.5){\footnotesize$\longleftarrow$}
  \textref h:C v:C \htext(109 41.5){\footnotesize$\longrightarrow$}

   \textref h:C v:C \htext(65 30){\footnotesize$ n-m-1$}
 \textref h:C v:C \htext(69 37){\footnotesize$\uparrow$}
  \textref h:C v:C \htext(69 23){\footnotesize$\downarrow$}

\linewd 0.1 \setgray 0
%lambda
\move(12 60)\lvec(8 60)
\move(12 100)\lvec(8 100)
\move(12 100)\lvec(12 104)
\move(48 100)\lvec(48 104)
%mu
\move(72 80)\lvec(68 80)
\move(72 100)\lvec(68 100)
\move(72 100)\lvec(72 104)
\move(112 100)\lvec(112 104)
%tau'
\move(30 54)\lvec(30 58)
\move(90 54)\lvec(90 58)
\move(30 50)\lvec(30 46)
\move(90 50)\lvec(90 46)
\move(70 50)\lvec(70 46)
%lambda'
\move(4 0)\lvec(8 0)
\move(4 40)\lvec(8 40)
\move(8 40)\lvec(8 44)
\move(48 40)\lvec(48 44)
%mu'
\move(112 44)\lvec(112 40)
\move(72 44)\lvec(72 40)
\move(72 40)\lvec(68 40)
\move(72 20)\lvec(68 20)

\textref h:C v:C \htext(73 112){\small$\tau$}
 \textref h:C v:C \htext(30 72){\small$\lambda$}
  \textref h:C v:C \htext(92 72){\small$\mu$}

\textref h:C v:C \htext(93 52){\small$\tau'$}
 \textref h:C v:C \htext(30 12){\small$\lambda'$}
  \textref h:C v:C \htext(92 12){\small$\mu'$} }

  \caption{The bijection $\varphi_2$ in Case $2$.}\label{fig-9-case2}
\end{figure}

Case $2$. For any $(\tau,\lambda,\mu )\in U_{n,m,k}$, we add a column
equals $2k-2$ to  $\lambda$ to obtain a new partition
 $\lambda'$ ,where $\ell(\lambda')=2k-2$ and $\lambda'_1 \leq
 m-2k+2$. By adding $2m-2k+2$ to $\tau$ we get
 $\tau'=((m-k+1)^2+k^2+m-k+1)$. And we remove a $2m$ part from $\mu$ and
 getting $\mu'$ whose largest part is less than or equal to $2m$ and length equals to
 $n-m-1$. This leads to the bijection $\varphi_2:U_{n,m,k}\rightarrow
 H_{n,m+1,k}$ as given by $(\tau,\lambda,\mu )\mapsto (\tau',\lambda',\mu')$.
 This case is illustrated in Figure \ref{fig-9-case2}.

\begin{figure}[h,t]
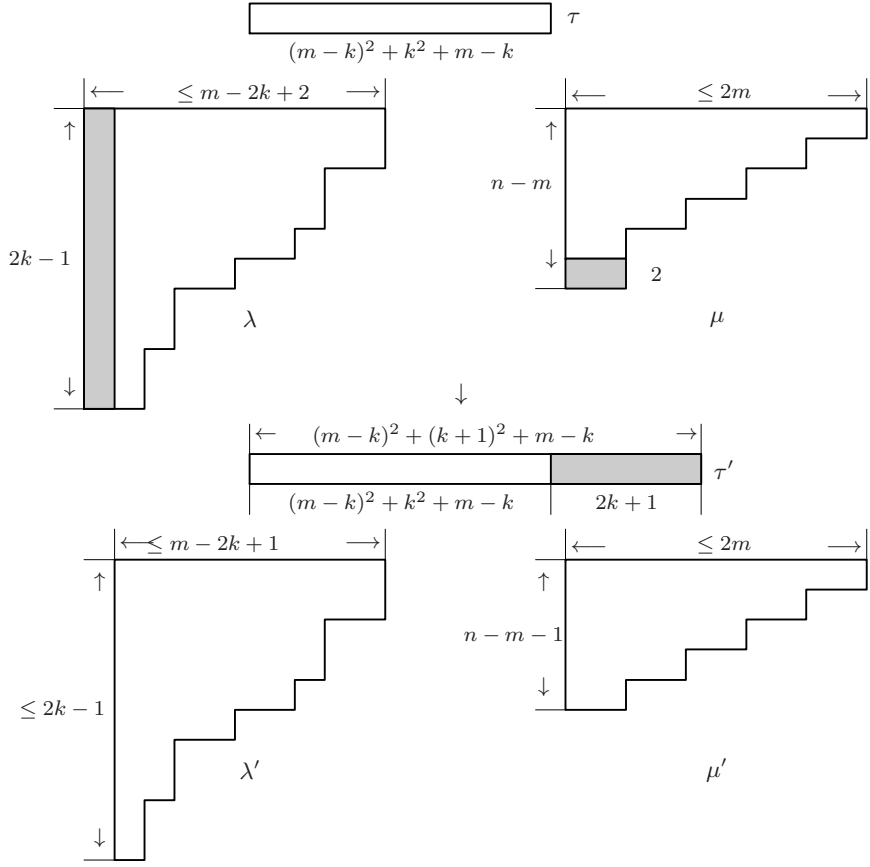

\centertexdraw{ \drawdim mm \linewd 0.25 \setgray 0
% \tau
\move(30 110)\lvec(70 110)\lvec(70 114)\lvec(30 114)\lvec(30 110)

\textref h:C v:C \htext(50 107.5){\footnotesize$(m-k)^2+k^2+m-k$}

% \lambda
\move(8 60)\lvec(16 60)\lvec(16 68)\lvec(20 68)\lvec(20 76)\lvec(28
76)\lvec(28 80)\lvec(36 80)\lvec(36 84)\lvec(40 84)\lvec(40 92)
\lvec(48 92)\lvec(48 100)\lvec(8 100)\lvec(8 60)
\move(8 60)\lvec(12 60)\lvec(12 100)\lvec(8 100)\lvec(8 60)\ifill f:0.8
\move(8 60)\lvec(12 60)\lvec(12 100)\lvec(8 100)\lvec(8 60)

 \textref h:C v:C \htext(29 102){\footnotesize$\leq m-2k+2$}
 \textref h:C v:C \htext(11 102){\footnotesize$\longleftarrow$}
 \textref h:C v:C \htext(45 102){\footnotesize$\longrightarrow$}

  \textref h:C v:C \htext(2 80){\footnotesize$ 2k-1$}
 \textref h:C v:C \htext(6 97){\footnotesize$\uparrow$}
 \textref h:C v:C \htext(6 63){\footnotesize$\downarrow$}

%\mu
\move(72 100)\lvec(112 100)\lvec(112 96)\lvec(104 96)\lvec(104
92)\lvec(96 92)\lvec(96 88)\lvec(88 88)\lvec(88 84)\lvec(80 84)\lvec(80 80)
\lvec(72 80)\lvec(72 100)
\move(72 80)\lvec(80 80)\lvec(80 76)\lvec(72 76)\lvec(72 80)\ifill f:0.8
\move(72 80)\lvec(80 80)\lvec(80 76)\lvec(72 76)\lvec(72 80)

\textref h:C v:C \htext(84 78){\footnotesize$2$}
 \textref h:C v:C \htext(93 102){\footnotesize$ \leq 2m$}
 \textref h:C v:C \htext(75 101.5){\footnotesize$\longleftarrow$}
  \textref h:C v:C \htext(109 101.5){\footnotesize$\longrightarrow$}

   \textref h:C v:C \htext(66 90){\footnotesize$ n-m$}
 \textref h:C v:C \htext(70 97){\footnotesize$\uparrow$}
  \textref h:C v:C \htext(70 80){\footnotesize$\downarrow$}

%downarrow
\textref h:C v:C \htext(58 62){\small$\downarrow$}

% \tau'
\move(30 50)\lvec(90 50)\lvec(90 54)\lvec(30 54)\lvec(30 50)
\move(70 50)\lvec(90 50)\lvec(90 54)\lvec(70 54)\lvec(70 50)\ifill f:0.8
\move(70 50)\lvec(90 50)\lvec(90 54)\lvec(70 54)\lvec(70 50)

\textref h:C v:C \htext(57 56.5){\footnotesize$(m-k)^2+(k+1)^2+m-k$}
 \textref h:C v:C \htext(32 56.5){\footnotesize$\leftarrow$}
  \textref h:C v:C \htext(88 56.5){\footnotesize$\rightarrow$}

\textref h:C v:C \htext(50 47.5){\footnotesize$(m-k)^2+k^2+m-k$}
%\textref h:C v:C \htext(57 44.5){\footnotesize$(m-k+1)^2+k^2+m-k+1$}
\textref h:C v:C \htext(80 47.5){\footnotesize$2k+1$}

% \lambda'
\move(12 0)\lvec(16 0)\lvec(16 8)\lvec(20 8)\lvec(20 16)\lvec(28
16)\lvec(28 20)\lvec(36 20)\lvec(36 24)\lvec(40 24)\lvec(40 32)
\lvec(48 32)\lvec(48 40)\lvec(12 40)\lvec(12 0)
%\move(8 0)\lvec(12 0)\lvec(12 40)\lvec(8 40)\lvec(8 0)\ifill f:0.8
%\move(8 0)\lvec(12 0)\lvec(12 40)\lvec(8 40)\lvec(8 0)

 \textref h:C v:C \htext(25 42){\footnotesize$\leq m-2k+1$}
 \textref h:C v:C \htext(15 42){\footnotesize$\longleftarrow$}
 \textref h:C v:C \htext(45 42){\footnotesize$\longrightarrow$}

  \textref h:C v:C \htext(5 20){\footnotesize$ \leq 2k-1$}
 \textref h:C v:C \htext(10 37){\footnotesize$\uparrow$}
 \textref h:C v:C \htext(10 3){\footnotesize$\downarrow$}

%\mu'
\move(72 40)\lvec(112 40)\lvec(112 36)\lvec(104 36)\lvec(104
32)\lvec(96 32)\lvec(96 28)\lvec(88 28)\lvec(88 24)\lvec(80 24)\lvec(80 20)
\lvec(72 20)\lvec(72 40)
%\move(72 40)\lvec(64 40)\lvec(64 20)\lvec(72 20)\ifill f:0.8
%\move(72 40)\lvec(64 40)\lvec(64 20)\lvec(72 20)

 \textref h:C v:C \htext(93 42){\footnotesize$\leq 2m$}
 \textref h:C v:C \htext(75 41.5){\footnotesize$\longleftarrow$}
  \textref h:C v:C \htext(109 41.5){\footnotesize$\longrightarrow$}

   \textref h:C v:C \htext(65 30){\footnotesize$ n-m-1$}
 \textref h:C v:C \htext(69 37){\footnotesize$\uparrow$}
  \textref h:C v:C \htext(69 23){\footnotesize$\downarrow$}

\linewd 0.1 \setgray 0
%lambda
\move(8 60)\lvec(4 60)
\move(8 100)\lvec(4 100)
\move(8 100)\lvec(8 104)
\move(48 100)\lvec(48 104)
%mu
\move(72 76)\lvec(68 76)
\move(72 100)\lvec(68 100)
\move(72 100)\lvec(72 104)
\move(112 100)\lvec(112 104)
%tau'
\move(30 54)\lvec(30 58)
\move(90 54)\lvec(90 58)
\move(30 50)\lvec(30 46)
\move(90 50)\lvec(90 46)
\move(70 50)\lvec(70 46)
%lambda'
\move(8 0)\lvec(12 0)
\move(8 40)\lvec(12 40)
\move(12 40)\lvec(12 44)
\move(48 40)\lvec(48 44)
%mu'
\move(112 44)\lvec(112 40)
\move(72 44)\lvec(72 40)
\move(72 40)\lvec(68 40)
\move(72 20)\lvec(68 20)

\textref h:C v:C \htext(73 112){\small$\tau$}
 \textref h:C v:C \htext(30 72){\small$\lambda$}
  \textref h:C v:C \htext(92 72){\small$\mu$}

\textref h:C v:C \htext(93 52){\small$\tau'$}
 \textref h:C v:C \htext(30 12){\small$\lambda'$}
  \textref h:C v:C \htext(92 12){\small$\mu'$} }

  \caption{The bijection $\varphi_3$ in Case $3$.}\label{fig-9-case3}
\end{figure}

Case $3$. For any $(\tau,\lambda,\mu )\in T_{n,m,k}$, we add $2k+1$ to
$\tau$ to get
 $\tau'=((m-k)^2+(k+1)^2+m-k)$. Then we remove a part equals $2$ from $\mu$ and
 getting $\mu'$, where $\ell(\mu')=n-m-1$. Similarly,by removing a part equals to $2k-1$
 from
 $\lambda$, we can get $\lambda'$, where $\ell(\lambda')\leq 2k-1$ and $\lambda'_1 \leq m-2k$.
 This leads to the bijection $\varphi_3:T_{n,m,k}\rightarrow
 K_{n,m+1,k+1}$ as given by $(\tau,\lambda,\mu )\mapsto
 (\tau',\lambda',\mu')$.
 Figure \ref{fig-9-case3} gives an illustration of the bijection $\varphi_3$.

The proof is completed by combining all the above bijections. \qed

Observe that the bijection $\varphi_1,\varphi_2$ and $\varphi_3$ preserve
the weight. The above theorem immediately leads to a recurrence
relation, which implies $q$-series identity \eqref{eq-q9a}.
To be more specific, we have the following corollary.

\begin{coro}
Let
\[
F_n(a,q)=\sum_{m,k\geq 0}(-1)^m\sum_{(\tau,\lambda,\mu )\in
P_{n,m,k}}a^n q^{\mid \tau \mid +\mid \lambda \mid+\mid \mu \mid},
\]
Then for any positive integer $n$, we have
\[
F_n(a,q)=-aq^{2n-1}F_{n-1}(a,q),\qquad n\ge2.
\]
Since $F_1(a,q)=-aq$, by iteration we find that
\[
F_n(a,q)=(-a)^n q^{n^2}.
\]
Summing over $n$, we arrive at identity \eqref{eq-q9a} of Andrews.
\end{coro}

\begin{figure}[h,t]
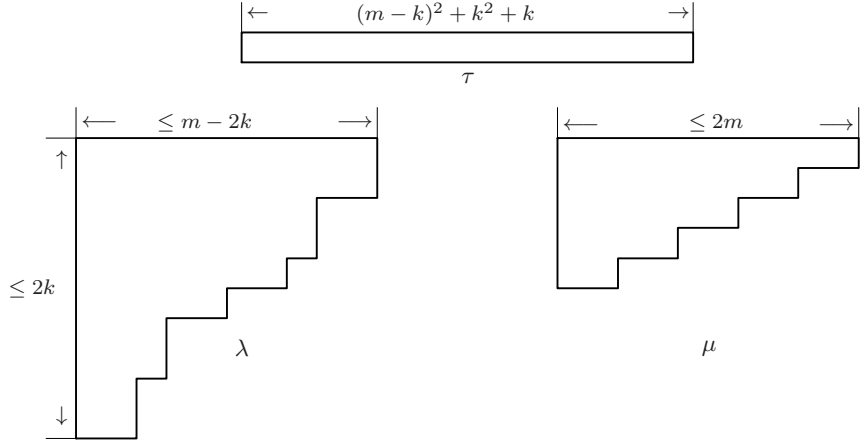

\centertexdraw{ \drawdim mm \linewd 0.25 \setgray 0

% \tau
\move(30 50)\lvec(90 50)\lvec(90 54)\lvec(30 54)\lvec(30 50)

\textref h:C v:C \htext(57 56.5){\footnotesize$(m-k)^2+k^2+k$}
 \textref h:C v:C \htext(32 56.5){\footnotesize$\leftarrow$}
  \textref h:C v:C \htext(88 56.5){\footnotesize$\rightarrow$}

% \lambda
\move(8 0)\lvec(16 0)\lvec(16 8)\lvec(20 8)\lvec(20 16)\lvec(28
16)\lvec(28 20)\lvec(36 20)\lvec(36 24)\lvec(40 24)\lvec(40 32)
\lvec(48 32)\lvec(48 40)\lvec(8 40)\lvec(8 0)

 \textref h:C v:C \htext(25 42){\footnotesize$\leq m-2k$}
 \textref h:C v:C \htext(11 42){\footnotesize$\longleftarrow$}
 \textref h:C v:C \htext(45 42){\footnotesize$\longrightarrow$}

  \textref h:C v:C \htext(2 20){\footnotesize$\leq 2k$}
 \textref h:C v:C \htext(6 37){\footnotesize$\uparrow$}
 \textref h:C v:C \htext(6 3){\footnotesize$\downarrow$}

%\mu
\move(72 40)\lvec(112 40)\lvec(112 36)\lvec(104 36)\lvec(104
32)\lvec(96 32)\lvec(96 28)\lvec(88 28)\lvec(88 24)\lvec(80 24)\lvec(80 20)
\lvec(72 20)\lvec(72 40)

 \textref h:C v:C \htext(93 42){\footnotesize$\leq 2m$}
 \textref h:C v:C \htext(75 41.5){\footnotesize$\longleftarrow$}
  \textref h:C v:C \htext(109 41.5){\footnotesize$\longrightarrow$}

\linewd 0.1 \setgray 0
%tau
\move(30 54)\lvec(30 58)
\move(90 54)\lvec(90 58)
%lambda
\move(4 0)\lvec(8 0)
\move(4 40)\lvec(8 40)
\move(8 40)\lvec(8 44)
\move(48 40)\lvec(48 44)
%mu
\move(112 44)\lvec(112 40)
\move(72 44)\lvec(72 40)

\textref h:C v:C \htext(60 48){\small$\tau$}
 \textref h:C v:C \htext(30 12){\small$\lambda$}
  \textref h:C v:C \htext(92 12){\small$\mu$} }
  \caption{Illustration of an element $(\tau,\lambda,\mu)\in
Q_{m,k}$.}\label{fig-10-case0}
\end{figure}

For another identity \eqref{eq-q10a} of Andrews, we give a sketch of the proof. Set
\[
G(q) =
\sum_{m,k\geq 0}\frac{(-a)^m
q^{(m-k)^2+k^2+k}}{(aq^2;q^2)_m}\left[\!
\begin{array}{c}
  m \\
  2k
\end{array}\!
\right]_q = \sum_{n\geq 0} (-a)^n q^{n^2},
\]
and the corresponding summand
\begin{equation}\label{eq-eq10-summation}
G_{m,k} = \frac{(-a)^m
q^{(m-k)^2+k^2+k}}{(aq^2;q^2)_m}\left[\!
\begin{array}{c}
  m \\
  2k
\end{array}\!
\right]_q = \sum_{n\geq 0} (-a)^n q^{n^2}.
\end{equation}
We also give a combinatorial interpretation of the summand $G_{m,k}$. 
Let $Q_{m,k}$ be a set of triple $(\tau, \lambda, \mu)$, where
$\tau$ is a partition with only one part equals 
$(m-k)^2+k^2+k$, $\lambda$ is a partition with no more than $2k$ parts and each part not exceeding $m-2k$,
and $\mu$ is a partition with only even parts not exceeding $2m$. 
Figure \ref{fig-10-case0} gives an illustration of an element of $Q_{m,k}$.
One can see that $G_{m,k}$ can be viewed as the weight of $Q_{m,k}$. 

Dividing $Q_{m,k}$ into a disjoint union of subsets
\[
Q_{n,m,k}=\{(\tau,\lambda,\mu )\in Q_{m,k}: \ell(\mu)=n-m\},
\]
where
$Q_{n,m,k}=\emptyset$ when $m=k=0$ or $2k>m>n$.  And let
\[
\begin{array}{l}
M_{n,m,k}=\{(\tau,\lambda,\mu )\in Q_{n,m,k}:\ell(\lambda)=2k,~m_2(\mu)=0\},\\[5pt]
S_{n,m,k}=\{(\tau,\lambda,\mu )\in Q_{n,m,k}:\ell(\lambda)=2k-1,~m_{2m}(\mu)=0\},\\[5pt]
L_{n,m,k}=\{(\tau,\lambda,\mu )\in Q_{n,m,k}:\ell(\lambda)\leq 2k-2,~m_{2m}(\mu)=0\}.
\end{array}
\]
We have the following combinatorial telescoping relation for $Q_{n,m,k}$.
\begin{theo}
For any positive integer $n$ and nonnegative integers $m$ and $k$, there is a bijection
\[
\begin{array}{l}
 \psi_{n,m,k}:  \\[5pt]
\quad Q_{n,m,k}\cup \{2n-1\} \times
Q_{n-1,m,k}\cup L_{n,m+1,k} \rightarrow \\[5pt]
  \quad M_{n,m,k}\cup M_{n,m+1,k}\cup  S_{n,m,k}\cup S_{n,m+1,k}\cup   L_{n,m,k}\cup L_{n,m+1,k+1}  \cup L_{n,m+1,k}.
\label{rec-q9}
\end{array}
\]
\end{theo}
\pf The proof of this theorem is similar to the proof of Theorem \ref{thm-mct}.\qed

The above theorem immediately leads to a recurrence relation as follows.

\begin{coro}
Let
\[
G_n(a,q)=\sum_{m,k\geq 0}(-1)^m\sum_{(\tau,\lambda,\mu )\in
Q_{n,m,k}}a^n q^{\mid \tau \mid +\mid \lambda \mid+\mid \mu \mid}.
\]
Then for any positive integer $n$, we have
\[
G_n(a,q)=-aq^{2n-1}G_{n-1}(a,q),\qquad n\ge2.
\]
\end{coro}

Since $G_1(a,q)=-aq$ and $G_0(a,q)=1$, by iteration we find that
\[
G_n(a,q)=(-a)^n q^{n^2}.
\]
Summing over $n$, we arrive at identity \eqref{eq-q10a} of Andrews.

%--------------------------------------------------------------------
\vspace{.2cm} \noindent{\bf Acknowledgments.}
We wish to thank Professor William Y.C. Chen for helpful comments and discussions. 
%This work was
%supported by the PCSIRT Project of the Ministry of Education, and the National Science Foundation of China.
%--------------------------------------------------------------------

\end{document}